\documentclass[preprint,11pt]{elsarticle}

\newcommand{\enumeration}{{\tt all-path}}
\newcommand{\BLM}{{\tt BLM}}
\newcommand{\basemodel}{{\tt base-model}}

\newcommand{\bap}{branch-and-price}
\newcommand{\bapC}{branch-and-cut-and-price}

\usepackage{amssymb}

\usepackage{amsmath}

\usepackage{optidef}

\usepackage{tikz}
\usetikzlibrary{arrows.meta}

\usepackage{float}

\usepackage{listings}

\usepackage{mathtools}

\usepackage{lineno}
\usepackage{fullpage}
\usepackage{palatino}
\usepackage{hyperref}
\usepackage{graphicx}

\usepackage{bm}

\def\equationrefname~#1\null{Equation (#1)\null}

\usepackage{natbib}
\usepackage[linesnumbered,ruled,vlined]{algorithm2e}

\SetCommentSty{mycommfont}
\DontPrintSemicolon
\SetKwIF{If}{ElseIf}{Else}{if}{}{else if}{else}{end if}
\SetKwFor{For}{for}{}{end for}%
\SetKwFor{ForEach}{foreach}{}{end foreach}%

\newcommand{\mv}{\mathcal{V}}
\newcommand{\ma}{\mathcal{A}}
\newcommand{\mk}{\mathcal{K}}
\newcommand{\mm}{\mathcal{M}}

\newproof{pf}{Proof}

\allowdisplaybreaks

\newtheorem{observation}{Observation} 

\begin{document}
\begin{sloppypar}

\title{Network Design with Service Requirements: Scaling-up the Size of Solvable Problems}


\author{Naga V. C. Gudapati}
\ead{chaitanya.gudapati@unibo.it}
\author{Enrico Malaguti}
\ead{enrico.malaguti@unibo.it}
\author{Michele Monaci}
\ead{michele.monaci@unibo.it}

\address{University of Bologna, BO, Italy, 40126}

\begin{abstract}
Network design, a cornerstone of mathematical optimization, is about defining the main characteristics of a network satisfying requirements on connectivity, capacity, and level-of-service. It finds applications in logistics and transportation, telecommunications, data sharing, energy distribution, and distributed computing.
In multi-commodity network design, one is required
to design a network minimizing the installation cost of its arcs and the operational cost to serve a set of point-to-point connections. The definition of this prototypical problem was recently enriched by additional constraints imposing that each origin-destination of a connection is served by a single path
satisfying one or more level-of-service requirements, thus defining the \textit{Network Design with Service Requirements} 
[Balakrishnan, Li, and Mirchandani. {\em Operations Research}, 2017].
These constraints are crucial, e.g., in telecommunications and computer networks, in order to ensure reliable and low-latency communication.
In this paper we provide a new formulation for the problem, where variables are associated with paths satisfying the end-to-end service requirements. We present a fast algorithm for enumerating all the exponentially-many feasible paths and, when this is not viable, we provide a column generation scheme that is embedded into a branch-and-cut-and-price algorithm. Extensive computational experiments on a large set of instances show that our approach is able to move a step further in the solution of the
Network Design with Service Requirements, compared with the {current state-of-the-art.}
\end{abstract}
 
\maketitle


\section{Introduction}
\label{Sec:1}

Network design is a cornerstone of mathematical optimization, as witnessed by the large amount of literature
on this topic. Indeed, historically it finds applications in logistics and transportation of goods and persons (\cite{MagnantiWong1984})
and, more recently, in telecommunications, data sharing, energy distribution, and distributed computing (\cite{GCF99}).

Network design is about defining the main characteristics of a network satisfying requirements on connectivity, capacity, and level-of-service.
Setting up the network induces some installation cost, while additional costs are incurred
when operating the service.  
It is quite common that a larger cost in the first term yields to a reduction in the latter, and
vice-versa.
Thus, the problem requires to find an equilibrium in the trade-offs between 
the installation and the operational costs.

A prototypical network design problem is the  multi-commodity network design, in which one is required
to design a network minimizing the installation cost of its arcs and the operational cost to serve a set of point-to-point connections, denoted as commodities. 
The solutions to this problem, however, can results in networks for which some commodities 
experience a low-quality connection with respect to some metric, e.g., distance or number of intermediate 
network nodes (hops) between origin and destination.
In some applications, this is a critical issue: for example in telecommunications, a common requirement 
consists of limiting the number of hops between origin and destination of any connection, as this has a direct effect on the latency of the communication. Similarly, in transportation networks, it is common to limit the distance traveled between origin-destination pairs, in particular when dealing with a public transport service
or when transporting perishable goods.

{Recently, \cite{BLM17} filled} this gap and introduced the \textit{Network Design with Service Requirements} (NDSR), {a network design
problem in which} additional constraints impose that each origin-destination is served by a single path
satisfying one or more level-of-service requirements. 
More specifically, each path must satisfy a maximum length with respect to a number of specified metrics. 
The problem asks to select some arcs to include in the network and to define, for each commodity, 
a path on the selected arcs and taking into account the mentioned level-of-service requirements. 
composed of {The objective is to minimize a cost function
consists of minimizing the total installation cost of the network arcs and of} the 
operational cost of the selected paths.
In that paper, the authors show that a model based on arc-flow variables can be hard to solve even for moderate-sized networks. Hence, through a wide polyhedral analysis they derive several families of valid inequalities, which can be exploited to strengthen the formulation. 
{The resulting model, combined with an 
effectiv heuristic algorithm, allows} to tackle larger 
instances of the problem.

In this manuscript, we propose a new model where variables are associated with paths satisfying the end-to-end service requirements. This way, many of the weaknesses of the arc-flow formulation are naturally 
overtaken without the need to recur to {cut separation} techniques. This desirable property comes at the cost of a formulation which is much larger, involving an exponential number of path variables. 
However, we show that for all the instances considered by \cite{BLM17}, we are indeed capable of quickly enumerating all the variables of the new formulation, thanks to an effective labelling algorithm, and 
to solve to proven optimality a much larger set of instances using a general-purpose ILP solver. 
In particular, our approach allows to solve a relevant fraction of {the} large instances
introduced by \cite{BLM17}, and to compute near-optimal solutions in the remaining cases, 
showing that the algorithm scales efficiently to larger size of the network. {In addition, we provide a new set of instances for which enumerating all the paths is not viable; for solving these large instances, we present a column generation scheme that is embedded into a full branch-and-cut-and-price algorithm.}

The paper is organized as follows. 
In the remainder of this section, we review some literature related to the problem at hand.
Section \ref{sec:description} formally describes
the problem, reviews a mathematical formulation from the literature, introduces a novel formulation, and compares the two models.
Section \ref{sec:branch-and-price} presents a solution approach based on branch-and-cut-and-price, describing column generation and the addition of 
valid inequalities. 
Section \ref{sec:computational} computationally
compares the performances of the proposed algorithms
with state-of-the-art approach on test instances from the literature.
Finally, in Section \ref{sec:conclusions} we present some conclusions.

\paragraph{Literature Review:}
There is a wide literature on network design problems, and many surveys have been published on these topics, see, e.g., 
\cite{MagnantiWong1984}, \cite{CRAINIC2000272}, and \cite{Wieberneit2007}. 
Depending on the specific application, different variants of these problems were considered. A notable field of research
involves the design of reliable and survivable networks, that has become a major objective for telecommunication
operators (see, \cite{KM05}).
In this context, one is required to define a {robust} network preserving a given 
connectivity level under possible failure of certain network components. 
There exist several ways to express the network robustness.
Under a stochastic paradigm, the network is required to remain operative either with a large 
probability (\cite{Song2013}, \cite{BARRERA2015132}) or after some recourse action has been implemented (\cite{LJUBIC2017333}).
Alternatively, more conservative approaches, imposing explicit redundancy in the definition of the network, have been 
considered in the literature; typically, one is required to design a network having two (edge) disjoint paths for each commodity
(\cite{MR05}, \cite{DBLP:journals/informs/AndreasS08}, \cite{DBLP:journals/jgo/AndreasSK08},
and \cite{10.1287/opre.1080.0579::Balakrishnan2009}), while \cite{GMS95} considered the case in which higher connectivity
requirements are imposed. 

Another class of related problems arises in applications where explicit constraints are imposed on the characteristics
of each path. A common requirement to guarantee the required quality of service is to limit the number
of hops of each path; this problem has been introduced by \cite{BA92}, while 
\cite{G98} presented a strong flow formulation that has been later adopted for many hop-constrained network
design problems.
In some cases, the resulting network is required to have a special structure (typically, a tree),
or survivability considerations have been added to the problem definition; see, e.g., \cite{BFGP13} and \cite{GLL15}.

{Our problem is closely related to the class of multi-commodity
flow problems (\cite{K78}) in which the network is given and commodities
compete for the use of the arcs, which have a limited capacity.
A branch-and-cut-and-price approach using path variables has been
proposed by \cite{BCV00}.}
Another relevant special case of NDSR arises when network design has to be defined for a unique commodity, and a single
metric has to be considered. The resulting budget constrained shortest path problem, introduced by \cite{J66}, 
is an NP-hard problem, and turns out to be a simplified version of a subproblem that we have to solve for generating columns, 
which takes {more than one metric} into consideration.

Finally, on the applications side, end-to-end service requirements have been 
considered by \cite{Barnhart1996}, \cite{Kim1999}, and \cite{Armacost2002},  where express 
delivery of parcels is optimized. Though service time is a key aspect in these applications, the special structure of the networks allows to avoid to explicitly 
impose these constraints.

\section{Problem description and formulation} \label{sec:description}

We now give a formal definition of the problem addressed in this
paper. 
We are given a directed graph $G = (\mv, \ma)$ where $\mv$ 
is the node set and $\ma$ is 
the arc set, and a set $\mk$
of commodities.
Each commodity $k \in \mk$ has associated a source node $s^k$ and a
sink node $t^k$. 
For each arc $a \in \ma$ there is an
activation cost $F_a$; in addition, using an arc $a$
for a commodity $k$ induces a flow cost $c_a^k$.
The problem asks to send, for each commodity $k$, 
one unit of flow on a single path $p^k$ 
from the source to the sink, 
by determining a set of arcs and the routing of the 
flows so that the sum of the 
activation and flow costs is a minimum.
In addition, there is a set $\mathcal{M}$ of metrics, 
that determines the feasibility of the path associated
with a given commodity $k$: for each metric $m$, we 
denote by $w_a^{km}$ the weight of arc $a$ with respect to 
the metric, and require that the sum of the weights
on arcs in $p^k$ does not 
exceed a given upper limit $W^{km}$.
We denote by ${\bf w}^k_a$ and ${\bf W}^k$ the corresponding $m$-dimensional vectors.

Throughout the paper, we assume that the 
graph includes no multiple arcs. This assumption is without loss
of generality, as multiple arcs with different costs or service
consumption for a given pair of nodes
can be handled by the addition of dummy nodes.
In addition, we assume that, 
for each commodity, at least
one feasible path exists, since otherwise the problem is 
clearly infeasible.

The problem reduces to the budget-constrained shortest path when
there is a single commodity and a single metric.
This shows that the problem is NP-hard.

The next section reports a descriptive formulation that
has been proposed in the literature, whereas
Section \ref{sec:path-based} introduces a novel formulation that will
be used in our solution scheme.

\subsection{Arc-flow formulation}\label{sec:arc-flow}

The following formulation has been proposed by
\cite{BLM17} and makes use of 
{activation} variables and flow variables. All variables
are binary and have the following meaning:

\begin{linenomath*}
\begin{equation*}
    z_{a}=
    \begin{cases}
     1 & \text{if}\ \text{arc}\ a \ \text{is selected} \\
     0 & \text{otherwise}
    \end{cases}
    \hspace*{12ex} (a \in \ma)
\end{equation*}
\end{linenomath*}
\begin{linenomath*}
\begin{equation*}
    y_{a}^{k}=
    \begin{cases}
     1 & \text{if}\ \text{commodity } k \text{ is routed on } \text{arc}\ a \\
     0 & \text{otherwise}
    \end{cases}
    \quad (a \in \ma, k \in \mk)
\end{equation*}
\end{linenomath*}

Then, the NDSR can be modelled using the following Integer
Linear Programming (ILP) formulation:

\begin{subequations}
    \begin{alignat}{2}
    & \min &\qquad \sum_{a \in \mathcal{A}}F_{a}z_{a} + \sum_{k \in \mathcal{K}}\sum_{a \in \mathcal{A}}c^{k}_{a}y^{k}_{a} \label{eq:baseNDSRobj}\\
    &\text{subject to} & \sum_{a \in \delta^{+}(v)} y_{a}^{k} - 
    \sum_{a \in \delta^{-}(v)}y_{a}^{k} &= 
    	\left\{ \begin{array}{ll}
	                    +1 & \hspace*{0.2cm} \mbox{$v = s^k$}\\
	                    -1 & \hspace*{0.2cm} \mbox{$v = t^k$ }\\
	                     0 & \hspace*{0.2cm} \mbox{$v \in \mv \setminus \{s^k,t^k\}$}
	                   \end{array} \right. 
	                   \quad k \in \mk 
    \label{eq:baseNDSRconstraint1}\\
&         &    \sum_{a \in \ma} w_{a}^{km}y_{a}^{k}                &\le W^{km} \hspace*{18.2ex} k \in \mk, m \in \mm \label{eq:baseNDSRconstraint5}\\
    &         &    y^{k}_{a} & \le z_{a} \hspace*{20.6ex} a \in \ma, k \in \mk \label{eq:baseNDSRconstraint4}\\
    &         &    z_{a}                               &\in \{0,1\} \hspace*{17.8ex} a \in \ma \label{eq:baseNDSRconstraint6A}\\
    &         &    y^{k}_{a}                               &\in \{0,1\} \hspace*{17.8ex} a \in \ma, k \in \mk. \label{eq:baseNDSRconstraint6Y}
    \end{alignat}
\end{subequations}

The objective function minimizes the sum of the 
{activation} and flow costs. 
Constraints \eqref{eq:baseNDSRconstraint1} impose flow conservation for each commodity
and node, whereas \eqref{eq:baseNDSRconstraint5} concern feasibility of the paths
with respect to the metrics, and inequalities \eqref{eq:baseNDSRconstraint4} 
force the activation of arcs that are used for routing a positive flow.
Finally \eqref{eq:baseNDSRconstraint6A} and \eqref{eq:baseNDSRconstraint6Y} define
the domain of the variables.
The arc-flow formulation has a polynomial size, as it includes $(|\mk|+1) \, |\ma|$ variables and
$|\mk| \, (|\mv| + |\ma| + |\mm|)$ constraints.

\subsection{Path-based formulation}\label{sec:path-based}

The novel ILP formulation that we propose includes the same binary {activation} variables of model \eqref{eq:baseNDSRobj}--\eqref{eq:baseNDSRconstraint6Y}, that select the arcs to be activated, {whereas
flow variables are replaced by path variables that are
defined as follows. Let $\mathcal{P}^{k}$ be} the set of all feasible paths for commodity $k$. 
For each commodity $k$ and each path $p \in \mathcal{P}^{k}$, 
let us introduce a binary path variable $x_p$ with the following meaning:
\begin{linenomath*}
    \begin{equation*}
        x_{p}= 
            \begin{cases}
                1, & \text{if }  \text{commodity } k\text{ is routed along path } p\\
                0 & \text{otherwise}
            \end{cases}\quad (k \in \mk, p \in \mathcal{P}^{k})
\end{equation*}
\end{linenomath*}
Let $c_{p}$ be the {flow} cost of the path $p$ for commodity $k$, defined as the sum of the flow costs of all the arcs in $p$. 
The problem can thus be modelled as follows:

\begin{subequations}
    \begin{alignat}{2}
    & \min &\qquad \sum_{a \in \mathcal{A}}F_{a}z_{a} + 
    \sum_{k \in \mathcal{K}} \sum_{p \in \mathcal{P}^{k}} c_p x_p \label{pathNDSRobj}\\
    &\text{subject to} & z_a - \sum_{p \in \mathcal{P}^{k}:\, a\in p} x_p &
    \ge 0 \hspace*{11.0ex} a \in \ma, k \in \mk \label{pathNDSRcons1}\\
    &         &    \sum_{p \in \mathcal{P}^{k}}x_{p} & = 1 \hspace*{11.0ex}  k \in \mk \label{pathNDSRcons2}\\
    &         &    z_{a}                               &\in \{0,1\} \hspace*{7.0ex} a \in \ma \label{pathNDSRcons3binZ}\\
    &         &    x_p                               &\in \{0,1\} \hspace*{7.0ex} p \in \mathcal{P}^{k}, k \in \mk. \label{pathNDSRcons3binX}
    \end{alignat}
\end{subequations}

The objective function minimizes {activation costs and flow} costs, which are here expressed in terms of path variables. Constraints \eqref{pathNDSRcons1} are the counterpart of \eqref{eq:baseNDSRconstraint4}, enforcing activation of arcs that are used by a path.  Constraints \eqref{pathNDSRcons2} ensure that, for every commodity, one feasible path is selected. Finally, \eqref{pathNDSRcons3binZ} and \eqref{pathNDSRcons3binX} define the domain of the variables.


\begin{observation}\label{obs:integrality}
The model obtained by relaxing integrality requirement \eqref{pathNDSRcons3binX}
admits an optimal integer solution.
\end{observation}
\begin{pf} 
Assume that an optimal solution for the relaxation is
given. 
{For a given choice of the $z$ variables, the} 
$x$ variables
associated with a commodity do not interact with those of a 
different commodity. Thus, we concentrate on a single commodity, 
say $k$, and assume that more than one path is selected for
that commodity, the sum of the values of the associated path variables being 1. 
By optimality of the initial solution, 
all the selected paths must have the same cost.
Hence, by increasing the value of one path variable to 1 and setting to 0 
all the remaining ones, we obtain a solution that has the same cost
as the original one. 
\hfill 
\end{pf}

\paragraph{Variable enumeration:}
We first observe that the 
path-based formulation has has $O(2^{|\mv|})$ variables
and $(|\ma|+1) \, |\mk|$ constraints, i.e., its size can be exponential
in the size of the instance. 
We now introduce an algorithm for enumerating all path variables; however, for large graphs, 
enumerating all paths can be challenging, and one may have to resort to column 
generation techniques, that will be discussed in Section \ref{sec:branch-and-price}.

Enumeration Algorithm \ref{algo:labelling-all} considers one commodity $k$ 
at a time and defines all simple paths from $s^k$ to $t^k$ that satisfy
resource constraints under all metrics.
The algorithm is inspired by the labelling method proposed by \cite{DB03} for the
{budget-constrained} shortest path problem.
In our algorithm, each label $\ell = \{u, c, {\bf w}\}$
represents a path from $s^k$ to $u$ having cost $c$ and using $w_m$ units of resources under
each metric $m$.
Each label is generated as unmarked, meaning that it has to be expanded, and then it is marked when  
considered for expansion.
Expansion of a label $\ell$ associated with a node $u$ consists in appending an arc $a = (u, v)$ 
to the current path. To this aim, we consider all the outgoing arcs from $u$ and, for each {neighbor node} 
 $v$ {not yet belonging to path $\ell$, we} check whether using the current label for reaching $v$ preserves feasibility with 
respect to the metrics. In this case, we define a new label 
$\ell' = \{v, c + c^{k}_a, {\bf w} + {\bf w}^k_a \}$, i.e., we
update the path cost and resource usage when using the current label for reaching $v$.
Eventually, node $v$ is inserted in set $T$, that includes all nodes 
associated with unmarked labels.
The algorithm terminates when $T = \emptyset$, meaning that no label can be further 
expanded, and returns all labels associated with node $t^k$.
Although a node can be inserted in and removed from $T$ more than once, the convergence of 
the algorithm is ensured by requiring simple paths, which is
checked in line \ref{test:feasibility}.

\begin{algorithm}[t!]
    \SetKwInOut{Input}{Input}
    \SetKwInOut{Output}{Output}
    \Input{$k$}
    $s := s^k, t := t^k, T := \{s\}, c := 0, {\bf w} := {\bf 0}$\;
    Define an unmarked label $\ell \coloneqq \{s, c, {\bf w}\}$ for node $s$\;
    \While{ $T \not= \emptyset$}
    {
        pick any $u \in T$\;
        $T := T \setminus \{u\}$\;
        \ForEach{{\em unmarked label} $\ell = \{u, c, {\bf w}\}$ {\em associated with node} $u$}
        {
            mark label $\ell$\;
            \ForEach{$a = (u, v) \in \delta^{+}(u)$} 
            { 
                \If{($v \notin$ {\em path} $\ell$) {\em and}
                $
                    ({\bf w} + {\bf w}^k_a \leq {\bf W}^k)
                $
                \label{test:feasibility}
                }
                {
                    define an unmarked label $\ell' = \{v, c + c^{k}_a, {\bf w} + {\bf w}^k_a \}$\;
                    \If{($v \not=t$)}
                    {
                        $T := T \cup\{v\}$  
                    }
                }
            }
        }   
    }
    \KwRet{{\em all labels associated with node} $t$}\;
    \caption{Compute all feasible paths for a fixed commodity}
    \label{algo:labelling-all}
\end{algorithm}

The above algorithm can be improved by pre-computing, for each metric $m \in \mm$, the shortest path 
from each node to $t^k$ when the cost of each arc $a$ is given by $w^{km}_a$. 
This figure can be used when checking feasibility of the new label in line \ref{test:feasibility}: by
adding this term to the left-hand-side of the inequality, we avoid generating labels that could not
be feasibly expanded to node $t^k$.

\subsection{Models comparison}
In this section we compare the two formulations in terms of their linear relaxations.
\begin{observation}\label{obs:equivalence}
Any feasible solution for the linear relaxation of the path-based formulation 
can be mapped to a feasible solution of the same cost of the linear relaxation of the
arc-based formulation, whereas the opposite does not hold.
\end{observation}

\begin{figure}[!ht]
\begin{center}
\begin{tikzpicture}
\begin{scope}[every node/.style={circle,thick,draw}]
    \node (s) at (0,0) {s};
    \node (1) at (4.0,2) {1};
    \node (t) at (8.0,0) {t};
\end{scope}
\begin{scope}[>={Stealth[black]},
              every node/.style={fill=none,circle},
              every edge/.style={draw=black,very thick}]
    \path [->] (s) edge[bend left=25] node[midway, above] {$[0,2]$} (1);
    \path [->] (s) edge[bend right=35] node[midway, below] {$[1,1]$} (t);
    \path [->] (1) edge[bend left=25] node[midway, above] {$[0,1]$} (t);
\end{scope}
\end{tikzpicture}
        \caption{Simple example for which the path-based formulation dominates the arc-flow formulation.}
        \label{fig:linear-relax}
\end{center}
\end{figure}
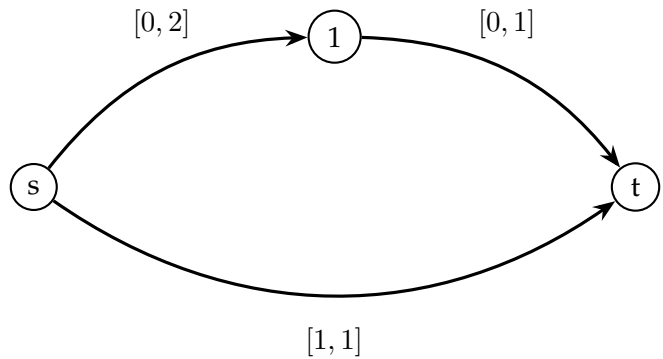

\begin{pf}
Let $z^*, x^*$ be a feasible solution of the linear relaxation of the
path-based formulation. We now define a solution $\widetilde z, \widetilde y$ that is
feasible for the linear relaxation of the arc-based formulation and has the same 
cost. First, we set $\widetilde z = z^*$. Then, for each arc $a \in \ma$ and
commodity $k \in \mk$, we set
$$
\displaystyle{\widetilde y_a^k = \sum_{p \in \mathcal{P}^{k}: a \in p} x^{*}_p}.
$$
It is straightforward to check that flow conservation constraints \eqref{eq:baseNDSRconstraint1} and feasibility requirements 
\eqref{eq:baseNDSRconstraint5} with respect to the metrics are satisfied
as $y$ variables are obtained as combination of feasible paths, whereas
constraints \eqref{eq:baseNDSRconstraint4} are implied by \eqref{pathNDSRcons2} and
by the definition of $\widetilde y_a^k$. The equivalence of the costs follows from
the definition of the cost of each path.

Figure \ref{fig:linear-relax} gives a small numerical example showing that the counterpart does not 
hold. The instance has no {flow} costs, a single commodity, and a single metric, for which the
capacity is $W = 2$. 
For each arc we report the {activation} cost and the weight with respect to the metric.
While there is a unique feasible path $p = \{(s, t)\}$ having cost 1, an optimal
solution to the linear relaxation of the arc-based formulation is given by 
$y_{s1} = y_{1t} = y_{st} = 1/2$ having cost 1/2.
\hfill
\end{pf}

The observation shows that the path-based formulation dominates the arc-based one
in terms of tightness of the associated linear relaxations.

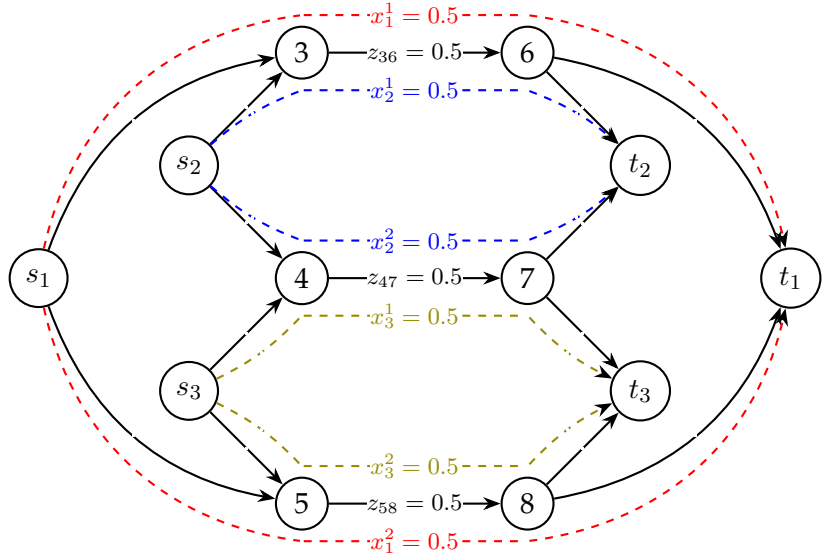
\begin{figure}[!ht]
        \begin{center}
                \begin{tikzpicture}
                        \begin{scope}[every node/.style={circle,thick,draw}]
                                \node (0) at (0.0,0.0) {$s_{1}$};
                                \node (1) at (2.0,1.5) {$s_{2}$};
                                \node (2) at (2.0,-1.5) {$s_{3}$};
                                \node (3) at (3.5,3.0) {3};
                                \node (4) at (3.5, 0.0) {4};
                                \node (5) at (3.5,-3.0) {5};
                                \node (6) at (6.5,3.0) {6};
                                \node (7) at (6.5,0.0) {7};
                                \node (8) at (6.5,-3.0) {8};
                                \node (9) at (8.0,1.5) {$t_{2}$};
                                \node (10) at (8.0,-1.5) {$t_{3}$};
                                \node (11) at (10,0) {$t_{1}$};
                        \end{scope}
        
                        \begin{scope}[every node/.style={inner sep=0pt,minimum size=1pt}]
                                \node (12) at (3.5,3.5) {};
                                \node (13) at (3.5,2.5) {};
                                \node (14) at (3.5,0.5) {};
                                \node (15) at (3.5,-0.5) {};
                                \node (16) at (3.5,-2.5) {};
                                \node (17) at (3.5,-3.5) {};
                                \node (18) at (6.5,3.5) {};
                                \node (19) at (6.5,2.5) {};
                                \node (20) at (6.5,0.5) {};
                                \node (21) at (6.5,-0.5) {};
                                \node (22) at (6.5,-2.5) {};
                                \node (23) at (6.5,-3.5) {};
                        \end{scope}
        
                        \begin{scope}[>={Stealth[black]},
                                      every node/.style={fill=white, inner sep=0pt,minimum size=1pt},
                                      every edge/.style={draw=black,thick}]
                            \path [->] (0) edge[bend left=30] node {} (3);
                            \path [->] (0) edge[bend right=30] node {} (5);
                            \path [->] (1) edge node {} (3);
                            \path [->] (1) edge node {} (4);
                            \path [->] (2) edge node {} (4);
                            \path [->] (2) edge node {} (5);
                           
                            \path [->] (6) edge[bend left=30] node {} (11); 
                            \path [->] (6) edge node {} (9); 
                            \path [->] (7) edge node {} (9); 
                            \path [->] (7) edge node {} (10); 
                            \path [->] (8) edge node {} (10); 
                            \path [->] (8) edge[bend right=30] node {} (11);
                            
                            \path [ thick, dashed] (0) edge[red, bend right=35] node {} (17); 
                            \path [ thick ,dashed] (17) edge[red] node {\footnotesize{$x_{1}^2 = 0.5$}} (23);
                            \path [thick ,dashed, ->] (23) edge[red, bend right=35] node {} (11);

                            \path [thick, dashed] (0) edge[red, bend left=35] node {} (12); 
                            \path [thick ,dashed] (12) edge[red] node {\footnotesize{$x_{1}^1 = 0.5$}} (18);
                            \path [thick ,dashed, ->] (18) edge[red, bend left=35] node {} (11);

                            \path [ thick, dashed] (1) edge[blue, bend right=10] node {} (14); 
                            \path [ thick ,dashed] (14) edge[blue] node {\footnotesize{$x_{2}^{2} = 0.5 $}} (20);
                            \path [thick ,dashed, ->] (20) edge[blue, bend right=10] node {} (9);

                            \path [ thick, dashed] (1) edge[blue, bend left=10] node {} (13); 
                            \path [ thick ,dashed] (13) edge[blue] node {\footnotesize{$x_{2}^{1} = 0.5 $}} (19);
                            \path [thick ,dashed, ->] (19) edge[blue, bend left=10] node {} (9);

                            \path [ thick, dashed] (2) edge[olive, bend right=10] node {} (15); 
                            \path [ thick ,dashed] (15) edge[olive] node {\footnotesize{$x_{3}^{1} = 0.5 $}} (21);
                            \path [thick ,dashed, ->] (21) edge[olive, bend right=10] node {} (10);

                            \path [ thick, dashed] (2) edge[olive, bend left=10] node {} (16); 
                            \path [ thick ,dashed] (16) edge[olive] node {\footnotesize{$x_{3}^{2} = 0.5$}} (22);
                            \path [thick ,dashed, ->] (22) edge[olive, bend left=10] node {} (10);

                            \path [->] (3) edge node {\footnotesize{$z_{36} = 0.5$}} (6);
                            \path [->] (4) edge node {\footnotesize{$z_{47} = 0.5$}} (7);
                            \path [->] (5) edge node {\footnotesize{$z_{58} = 0.5$}} (8); 
                            
                        \end{scope}
                        \end{tikzpicture}
        \end{center}
        \caption{Fractional solution of the linear relaxation of the path-based formulation.}
        \label{fig:fractional-LP}
\end{figure}

The structure of feasible solutions for the linear relaxation of the 
arc-based formulations was analyzed by \cite{BLM17}, showing that fractional
solutions may arise for two main reasons:
\begin{itemize}
    \item for a given commodity, the model may route part of the flow on a path
    that is less expensive but infeasible with respect to the metric
    requirements (see again Figure \ref{fig:linear-relax});
    \item arc {activation} variables can be set at a fractional value to allow sharing the
    {activation} cost of some arcs among different paths associated with 
    different commodities.
\end{itemize}
Accordingly, \cite{BLM17} introduced different families 
of valid inequalities to cut some of these solutions.
The first type of fractionalities do not appear in the path-based formulation,
in which feasibility of the paths is enforced when defining the variables; thus,
adding similar inequalities would be useless. On the other hand, the second type of fractionality 
may affect the path-based formulation as well, as shown in Figure \ref{fig:fractional-LP}.
In this example, there are three commodities, no flow costs and 
activation costs equal to one for arcs $(3,6), (4, 7), (5,8)$ and zero for 
the remaining arcs. The figure shows an optimal solution of the linear relaxation of
the path-based formulation, where the flow of each commodity is split into two paths,
the costly arcs are activated at value 0.5 and the resulting cost is 3/2. On the other hand,
any integer feasible solution has a cost at least equal to 2.
For this reason, in our approach we consider the possibility to add some classes of valid 
inequalities of the second type.

\section{Branch-and-cut-and-price approach}\label{sec:branch-and-price}

In this section, we introduce an exact algorithm based on the path-formulation that
can be used when enumerating all paths is unpractical.
The algorithm adopts a branch-and-bound strategy and solves, at each node, the linear relaxation 
of the model by means of column generation techniques.
The basic scheme is 
possibly enriched by the addition of valid inequalities, that do not change
the structure of the method, thus resulting in a robust branch-and-cut-and-price
algorithm.

\subsection{Column generation and labelling}

Column generation is an iterative scheme used for solving linear models with an
exponentially large number of variables. At each iteration, a 
{\em restricted master problem} including a subset of the variables is solved,
and its dual solution is used to determine new variables (if any) that have
to be added to the formulation in order to converge to an optimal solution.

In our setting we assume without loss of generality that constraints
\eqref{pathNDSRcons2} are rewritten as inequalities. 
At each iteration, the restricted master includes all the $z$ variables, 
and a non-empty subset $\widetilde{\mathcal{P}^{k}}\subseteq \mathcal{P}^{k}$ of
path variables for each commodity $k$ (notice that by construction
the restricted master always includes a feasible solution).
Assume that the restricted master has been solved to optimality, and let
$\gamma^{k}_{a}$ and $\rho^{k}$ be optimal non-negative dual variables associated with
constraints \eqref{pathNDSRcons1} and \eqref{pathNDSRcons2}, respectively.
The {\em reduced cost} for a path variable $x_p$ for a commodity $k$ is
given by 
$$
\overline c_p = 
c_p + \sum_{a \in p} \gamma^{k}_{a} - \rho^{k} = 
\sum_{a \in p} \big(c^k_a + \gamma^{k}_{a} \big) - \rho^{k} = 
\sum_{a \in p} \widetilde c^k_{a} - \rho^{k},
$$
where the arc costs are $\widetilde c^k_{a} = c^k_{a} + \gamma^{k}_{a}$.
Thus, the {\em pricing problem} for a given commodity $k$ is to find a feasible path 
whose reduced cost is negative, and can be formulated as a 
{\em budget-constrained shortest path problem} under costs 
$\widetilde c^k _{a}$ and resources defined by the metrics.
If the cost of this shortest path is strictly smaller than $ \rho^{k}$, the
corresponding path variable is added to the restricted master, and the process is iterated; if no 
path variable is generated for any commodity, the optimal solution of the current
restricted master is an optimal solution for the linear relaxation of the
problem.

\paragraph{Solution of the budget-constrained shortest path problem:}
Enumeration Algorithm \ref{algo:labelling-all} can be modified to compute 
the shortest path under resource constraints for a given commodity $k$, a problem which 
is NP-hard even if the graph is acyclic and $|\mm| = 1$ (see, \cite{GJ79}).
The resulting Algorithm \ref{algo:labelling-best} differs from the 
enumeration one starting from line \ref{test:dominance1}, where a dominance check 
aimed at avoiding expansion of suboptimal paths is introduced. More precisely, 
label $\ell'$ is dominated by another label $\ell''$ associated with the same node if
its cost and its resource usage are larger then or equal to the cost and usage of $\ell''$. In this case $\ell'$ is marked. Vice-versa, it may also happen that $\ell'$ dominates $\ell''$, in which case we 
mark $\ell''$.
Node $v$ is inserted in set $T$ only if label $\ell'$ remains unmarked.
The algorithm returns a unique path, corresponding to the label with minimum cost among all those associated with node $t^k$.

\begin{algorithm}[!ht]
    \SetKwInOut{Input}{Input}
    \SetKwInOut{Output}{Output}
    \Input{$k$}
    $s := s^k, t := t^k, T := \{s\}, c := 0, {\bf w} := {\bf 0}$\;
    Define an unmarked label $\ell \coloneqq \{s, c, {\bf w}\}$ for node $s$\;
    \While{ $T \not= \emptyset$}
    {
        pick any $u \in T$\;
        $T := T \setminus \{u\}$\;
        \ForEach{{\em unmarked label} $\ell = \{u, c, {\bf w}\}$ {\em associated with node} $u$}
        {
            mark label $\ell$\;
            \ForEach{$a = (u, v) \in \delta^{+}(u)$} 
            { 
                \If{($v \notin$ {\em path} $\ell$) {\em and}
                $
                    ({\bf w} + {\bf w}^k_a \leq {\bf W}^k)
                $
                }
                {
                    define an unmarked label $\ell' = \{v, c + \widetilde c^{k}_a, {\bf w} + {\bf w}^k_a \}$\;
                    \If{($\ell'$ {\em is dominated by a label} $\ell''$ {\em associated with node} $v$)
                    \label{test:dominance1}
                    }
                    {
                        mark label $\ell'$\;
                    }
                    \If{($\ell'$ {\em dominates a label} $\ell''$ {\em associated with node} $v$)}
                    {
                        mark label $\ell''$\;
                    }
                    \If{($v \not=t$) {\em and} ($\ell'$ {\em is unmarked})}
                    {
                        $T := T \cup\{v\}$  
                    }
                }
            }
        }   
    }
    \KwRet{{\em the unmarked label with minimum cost $c$ associated with node} $t$}\;
    \caption{Compute a constrained shortest path for a fixed commodity}
    \label{algo:labelling-best}
\end{algorithm}

\subsection{Branching scheme}

In our branching scheme we always select a $z$ variable for branching.
According to Observation \ref{obs:integrality}, at each node where all the $z$ variables attain
integer values, there exists an optimal solution in which all the $x$ variables are integer as well.
Notice that this is the solution returned by solving the restricted master problem by means of the
simplex algorithm.

A positive effect of this branching strategy is that it does not affect the structure 
{of the pricing subproblem.}
This is a crucial property for designing an effective branch-and-price algorithm, as it allows
to solve the column generation subproblem throughout all the branching tree by means of the
same effective labelling algorithm used at the root node.
Clearly, imposing $z_a = 1$ for some $a \in \ma$ has no direct effect in the {pricing}.
Conversely, {when imposing $z_a = 0$,
in the pricing subproblem we simply forbid the use of arc $a$ when generating new path variables,} which can be easily handled by setting $\ma = \ma \setminus \{a\}$.


\subsection{Adding valid inequalities}

In order to tighten the formulation and increase the dual bound at each node, we can add
valid inequalities that cut fractional solutions in which arc {activation} variables are set at a fractional 
value to allow sharing the {activation} cost of some arcs among different paths.

To this aim, we adapt to our model some of the inequalities introduced by \cite{BLM17} for
the arc-flow formulation.
These inequalities are obtained by analyzing the structure of the graph $G$ and by deriving relationships between pairs of arcs $(a,b)$ when routing the flow of a commodity $k$, namely:
\begin{itemize}
        \item 
\emph{OR} relationships, occurring when no more than one arc of pair $(a,b)$ can be used to route flow from $s^{k}$ to $t^{k}$;
       \item 
\emph{IF} relationships, occurring when the flow through arc $a$ must also be routed through $b$; and
        \item  
\emph{CUT} relationship, occurring when at least one between $a$ and $b$ must be used to route the flow. 
\end{itemize}
These relationships are then used to derive conditions that link the 
{activation} variable of an arc with the
flow variables associated with the same arc and different commodities.
By using the arc-flow variables, all these inequalities have the following general structure
\begin{align*}
         \sum_{(a,k) \in C} z_{a} - \sum_{(a,k) \in C} y^{k}_{a} \ge q,
\end{align*}
where $C$ is a set of arc-commodity pairs and $q$ is a scalar number.

By translating these conditions in terms of the path variables, we obtain
\begin{align}\label{eq:validcut}
         \sum_{(a,k) \in C} z_{a} - \sum_{(a,k) \in C}\sum_{p \in {\cal P}^k: a \in p}  x_{p} \ge  q  
\end{align}
which can be enforced in the path-based formulation.

As it happens for the branching conditions, the addition of the inequalities above does
not affect the structure of the pricing problem at a generic node of the branching tree.
Indeed, for a given commodity $k$, constraint \eqref{eq:validcut} only affects those paths 
that contain an arc $a$ such that pair $(a,k) \in C$. For each such path, the reduced cost 
of the associated variable is thus
$\overline c_p = 
\sum_{a \in p} \big(c^k_a + \gamma^{k}_{a} \big) + \phi^C - \rho^{k}
$
where $\phi^C$ is the dual variable associated with constraint \eqref{eq:validcut}.
More in general, given a collection $\cal{C}$ of inequalities, the reduced cost of
a path associated with commodity $k$ is

\begin{align*}
\overline c_p = 
\sum_{a \in p} \big(c^k_a + \gamma^{k}_{a} + \sum_{C \in {\cal C}: (a,k)\in C} \phi^C \big)- \rho^{k}
\end{align*}

Hence, the only effect of additional inequalities on the shortest path computation is on 
the definition of arc costs $\widetilde c^k_{a}$, which now include the dual variables of these 
constraints as well.

This allows us to solve the column generation subproblem with no modification of the 
labelling algorithm even after the addition of valid inequalities. The resulting
algorithm is then a robust branch-and-cut-and-price.

\section{Computational experiments}\label{sec:computational}

In our computational experiments we explore three directions. First, we compare the computational
performance of the path-based formulation with the arc-based formulation. Our second order of business is to 
determine the features of the instances for which full enumeration of all feasible paths
is possible, and when instead one has to resort to column generation. In this case, the solver
cannot be used as a black box, and the addition of valid inequalities may be an effective option
for accelerating the solution process. Finally, we evaluate the effect of adding valid inequalities to the path-based formulation, in terms of bound given by the linear relaxation and overall performance of the algorithm.

Unless specified, all algorithms were run on an
AMD Ryzen Threadripper 3960X running at 3.8 GHz in single-thread mode, with a time
limit of 1 hour per instance.
All algorithms were implemented in C++.
Both the arc-flow and the path-based formulations were solved
using Gurobi version 9.1.1 as ILP solver, whereas
the branch-and-cut-and-price was implemented on top of the SCIP optimization suite
(version 7.0.1 with its default {SoPlex} solver), which
allows to embed a column generation scheme within the enumeration process
(see \cite{GamrathEtal2020ZR}).

\subsection{Instances from the literature}\label{sec:lit-instances}

We now describe a benchmark of instances that has recently been introduced by \cite{BLM17},
who kindly provided us the code for generating the numerical data.
Each instance is characterized by the following parameters: the number of nodes $|\mv|$, 
number of arcs $|\ma|$, and number of commodities $|\mk|$. Nodes are randomly located on a 
rectangular grid and are connected by a spanning arborescence; then, $|\ma| - |\mv| + 1$ arcs
are added to the arc set, making sure that the resulting network contains one directed path
for each pair of nodes. 
The source and terminal node for each commodity are randomly selected in $\mv$.
The activation cost of each arc depends on the
euclidean distance between the two endpoints and on a random parameter. A parameter $\gamma$
governs the ratio between flow costs and activation costs.
Coefficients ${\bf w}^k_a$ for a given arc $a \in \ma$ are negatively correlated to the activation cost
$F_a$ through a parameter $\beta$ and a random term. 
All instances consider $|\mm|=2$ metrics. Weight limits for each commodity $k$ and each metric equal the length (using arc weights as lengths) of the $q$-th shortest  path from $s^k$ to $t^k$, where $q$ is a random parameter having uniform distribution in an interval of size $\Delta Q$ centered in $Q_{avg}$.
A particular combination of network size ($|\mv|$, $|\ma|$, and $|\mk|$),
cost structure and service requirements ($\beta$, $\gamma$, $Q_{avg}$, and 
$\Delta Q$) is referred to as a scenario. Overall, \cite{BLM17} 
defined 18 scenarios:
the first seven scenarios share the same default values of the parameters for cost
structure and service requirements, while considering varying network sizes
ranging from 30 nodes, 120 arcs, and 90 commodities to 
50 nodes, 250 arcs, and 150 commodities.
Scenarios 8-15 are all defined with a fixed network size
($|\mv| = 50$, $|\ma| = 200$, and $|\mk| = 150$) and different
cost structure and service requirements.
Finally, the last 3 scenarios have the default values of the 
parameters defining cost structure and service requirements and are 
characterized by larger size of the network, up to 
80 nodes, 320 arcs, and 240 commodities.
For each scenario, five instances were generated, for a total of  90 instances. 
Although the original set of instances is not available, we generated 18 scenarios for a total of 90 instances by using the same parameters used by \cite{BLM17}.
In Table \ref{tab:results-literature}, we will refer to each scenario as
$|\mv| / |\ma| / |\mk| / \beta \, \gamma \, Q_{avg} \, \Delta Q$
where the last four parameters take values in $\{L, M, H\}$ to denote 
low, medium and high figures, respectively.

\subsection{Results on the instances from the literature}\label{sec:results-lit-instances}

Table \ref{tab:results-literature} gives the results of
computational experiments on the 90 instances derived from the 18
scenarios described above; instances are grouped by scenario, i.e., every row
reports aggregate results for five instances.
The table compares the following approaches:
\begin{itemize}
    \item \basemodel\ corresponds to the direct application of general-purpose ILP 
    solver Gurobi to the arc-flow formulation;
    \item \BLM\ is the \textit{composite} algorithm proposed by \cite{BLM17}, and implements 
    a branch-and-cut scheme built on top of
    the general-purpose ILP solver
    Cplex 12.5.1 for solving the arc-flow formulation. 
    The algorithm includes separation of 
    several families of valid inequalities and 
    an effective LP-based heuristic algorithm that is executed at the root node.
    All these figures are taken from \cite{BLM17}, and correspond to experiments executed on a Intel core i5 using an integrality gap for
    early termination equal to 0.1\%;
    \item \enumeration\ denotes the algorithm obtained by enumerating 
    all feasible paths through Algorithm \ref{algo:labelling-all} and solving the resulting path-based
    formulation using the Gurobi ILP solver.
    This approach does not include cut separation nor column generation, 
    allowing us to use the solver as a black box, so as to exploit its full capabilities. 
\end{itemize}

For each solution approach, the table reports the number of instances
solved to proven optimality, the average percentage gap, and the average computing time (in seconds, with respect to instances that are solved to optimality only). 
For a given instance of the problem, let $L$ and ~$U$ be the best
lower and upper bound, respectively, produced by an algorithm; the resulting percentage gap is computed as 
$100 \frac{U-L}{U}$. 
For algorithm \BLM, detailed computational results are only 
available for the instances of the first 7 scenarios.
In addition observe that, for some scenarios, this algorithm solves all the associated
5 instances to optimality though returning a strictly positive
percentage gap, due to the tolerance value
that is used within the algorithm.
Finally, for algorithm \enumeration\ we also report the number of path variables
enumerated by the labelling algorithm. The enumeration time is always very small
(at most 0.5 seconds) and it is included in the computing time of the algorithm.

\begin{table}[h!]
\centering      
\small
\tabcolsep=2.7pt
\begin{tabular}{rr|rrr|rrr|rrrr}
          &                    & \multicolumn{3}{c|}{\basemodel} & \multicolumn{3}{c|}{\BLM} & \multicolumn{4}{c}{\enumeration} \\
          &                    & \multicolumn{3}{c|}{} & \multicolumn{3}{c|}{(\cite{BLM17})} & \multicolumn{4}{c}{} \\
\multicolumn{2}{c|}{scenario} & \# opt & \% gap &    time & \# opt & \% gap & time & \# opt & \% gap &   time & \# path \\
\hline
1         &     30/120/90/MMMM &      5 &   0.00 &  446.75 &      5 &   0.02 &  175 &      5 &   0.00 &   1.41 &    895 \\
2         &    40/160/120/MMMM &      4 &   0.21 & 1352.68 &      4 &   0.18 &  133 &      5 &   0.00 &   7.66 &   1098 \\
3         &    50/150/150/MMMM &      5 &   0.00 &  776.90 &      5 &   0.08 &  230 &      5 &   0.00 &   3.24 &   1409 \\
4         &    50/200/100/MMMM &      2 &   1.36 & 2953.33 &      4 &   0.43 & 1055 &      5 &   0.00 &  41.71 &    956 \\
5         &    50/200/150/MMMM &      1 &   5.93 & 2817.12 &      2 &   1.33 &  350 &      5 &   0.00 & 478.80 &   1455 \\
6         &    50/200/200/MMMM &      0 &   3.81 &      -- &      3 &   0.45 &  735 &      5 &   0.00 & 324.22 &   1898 \\
7         &    50/250/150/MMMM &      0 &   9.20 &      -- &      1 &   2.61 & 2631 &      4 &   0.38 & 171.85 &   1388 \\
8         &    50/200/150/LMMM &      2 &   3.09 & 2455.09 &        &   0.50 &      &      5 &   0.00 &  32.90 &   1249 \\
9         &    50/200/150/HMMM &      0 &  12.02 &      -- &        &   2.00 &      &      3 &   1.64 & 365.97 &   1795 \\
10        &    50/200/150/MLMM &      0 &   9.02 &      -- &        &   0.90 &      &      5 &   0.00 & 639.56 &   1455 \\
11        &    50/200/150/MHMM &      1 &   6.60 & 3557.67 &        &   0.10 &      &      5 &   0.00 & 782.21 &   1455 \\
12        &    50/200/150/MMLM &      5 &   0.00 & 1009.59 &        &   0.10 &      &      5 &   0.00 &   1.51 &    804 \\
13        &    50/200/150/MMHM &      0 &  10.68 &      -- &        &   1.90 &      &      4 &   0.64 & 345.37 &   2085 \\
14        &    50/200/150/MMML &      0 &   6.26 &      -- &        &   0.70 &      &      5 &   0.00 & 305.79 &   1421 \\
15        &    50/200/150/MMMH &      1 &   2.79 & 1902.54 &        &   0.70 &      &      5 &   0.00 & 100.74 &   1153 \\
16        &    60/240/180/MMMM &      0 &  10.68 &      -- &        &   0.70 &      &      5 &   0.00 & 603.16 &   1711 \\
17        &    70/280/210/MMMM &      0 &  11.91 &      -- &        &   2.50 &      &      2 &   0.54 & 686.22 &   1961 \\
18        &    80/320/240/MMMM &      0 &  15.86 &      -- &        &   2.40 &      &      2 &   1.25 & 951.94 &   2307 \\
\hline     
\multicolumn{2}{c}{summary}    &     26 &   6.08 & 1371.97 &     45$^*$ &   0.98 &      &     80 &   0.25 & 288.22 &   1472 \\ 
\end{tabular}
\caption{Results on instances from the literature.}\label{tab:results-literature}
\end{table}

The results confirm the outcome of the computational experiments reported by \cite{BLM17} for
the first seven scenarios, i.e., that algorithm \BLM\ outperforms the \basemodel, which can solve only 
small instances and has large percentage gaps for most unsolved scenarios.
Instead, results borrowed from \cite{BLM17} show that the addition of valid inequalities and the use of an effective heuristic
yields to an algorithm which is able to solve
24 instances out of 35, with average percentage gap 
equal to 0.73. Both these approaches are dominated by algorithm \enumeration, which
solves all but one instances in the first seven scenarios, and has a percentage
gap equal to 0.38 for the remaining instance.
This is due to the fact that the formulation is tight and 
that, for these instances, the number of path variables does
not grow up: this number is always smaller than 2000, which makes 
the model solvable with a limited 
computational effort.
All instances for scenarios 1 and 3 are solved
by both \BLM\ and \enumeration: for these scenarios, the average computing time
of the former is two orders of magnitude slower than
the latter (although \BLM\ was executed on a slightly
slower machine and used a different ILP solver).

For what concerns the instances in scenarios 8-15, the
performances of algorithm
\enumeration\ remain satisfactory. The algorithm solves
37 of the 40 associated instances, and has an average percentage gap equal to 0.28.
Finally, for very large instances (scenarios 16 to 18), 
the algorithm solves 9 instances out of 15 and
has an average percentage gap equal to 0.60.
Overall, our algorithm solves to proven optimality 
almost 90\% of the instances with an average percentage gap of
0.25.
\cite{BLM17} do not report detailed
results for all scenarios, but instead mention that \BLM\ only solves 
45 instances 
(for this reason this figure is marked with an asterisk in the 
summary line of the table) and 
has an average percentage gap of 0.98.

\subsection{Results on additional instances}\label{sec:results-new-instances}

The results in Table \ref{tab:results-literature} show that, for the instances from the literature,
the number of feasible paths is quite small. Thus, not surprisingly, the \enumeration\ approach
is always better than \basemodel\ and \BLM.
Our second set of experiments is aimed at evaluating the limits of applicability of explicit
enumeration of all path variables, and the alternative use of the branch-and-price 
algorithm described in Section \ref{sec:branch-and-price} when enumeration
is unpractical. 
Hence, we generated additional instances derived from the 
instances in scenarios 1--7, in which the number of feasible paths is increasing.
To minimize the number of parameters for defining the additional instances, we simply introduce a parameter $\alpha \geq 1$ that is used
to scale each upper limit $W^{km}$ for a commodity $k$
and metric $m$. This has the effect to make less binding the
constraints defining the feasibility of a path with respect to
the metrics.

\begin{table}[b!]
\centering      
\tabcolsep=3.0pt
\begin{tabular}{r|rrr|rrrr|rrrr}
         & \multicolumn{3}{c}{\basemodel} &        \multicolumn{4}{c}{\enumeration} &               \multicolumn{4}{c}{\bap} \\
$\alpha$ & \# opt & \% gap &         time & \# opt & \% gap &    time & \# path & \# opt & \% gap &  time  & \# path \\
\hline     
    1.00 &     17 &   2.93 &      1191.34 &     34 &   0.05 &  146.25 &    1300 &     31 &   0.38 & 401.20 &   1116 \\
    1.25 &      3 &   8.45 &      1680.47 &     21 &   1.45 &  410.74 &    6428 &     16 &   2.65 & 816.47 &   5896 \\
    1.50 &      6 &   6.01 &       976.83 &     20 &   1.72 &  514.40 &  33,178 &     14 &   2.64 & 667.09 & 10,816 \\
    1.75 &      9 &   3.97 &       903.57 &     19 &   1.71 &  480.65 & 169,286 &     17 &   2.15 & 539.52 & 11,209 \\
    2.00 &     14 &   2.40 &       798.33 &     17 &   1.97 &  449.87 & 855,441 &     21 &   1.55 & 830.85 & 10,110 \\
    2.25 &     18 &   1.53 &       613.12 &     14 &     -- &  673.01 &      -- &     23 &   1.19 & 522.94 &   9115 \\
    2.50 &     22 &   0.91 &       654.06 &      9 &     -- &  957.18 &      -- &     26 &   0.80 & 475.83 &   7743 \\
    2.75 &     23 &   0.70 &       554.46 &      5 &     -- &  998.54 &      -- &     27 &   0.68 & 577.48 &   7479 \\
    3.00 &     25 &   0.61 &       470.26 &      3 &     -- & 1627.80 &      -- &     28 &   0.58 & 628.28 &   7048 \\
\hline
\end{tabular}
\caption{Results on additional instances.}\label{tab:results-additional}
\end{table}

Table \ref{tab:results-additional} reports aggregated results, 
summarizing 35 instances per line, obtained
with different values of $\alpha$ ranging from $1.00$ to $3.00$.
We compare the \basemodel, the \enumeration\ approach, and the 
\bap\ algorithm and report, for each solution method, the number of optimal solutions,
the average percentage gap and the average computing time (with respect to
instances solved to optimality only). For \enumeration\ we also report
the total number of feasible paths; this figure is averaged over all the 
35 instances of a line, provided that enumeration of all paths
was completed within the time limit for all the instances.
Finally, for \bap\ we give the average number of path variables that have been generated
during the execution of the algorithm (with respect to instances solved
to optimality only).

The results in Table \ref{tab:results-additional} show that, for values of
$\alpha < 2$, the total number of paths is still manageable (below 200,000) 
and \enumeration\ remains the best option. 
Conversely, for larger values of $\alpha$, in many cases 
enumerating all path variables within the time limit is not possible or the path-based formulation has too many variables, and hence a method based
on column generation is advisable. Indeed, 
for $\alpha=2$, \bap\ solves
21 instances compared to the 17 solved 
by \enumeration, and this gap increases for larger values of $\alpha$.
Finally, we observe that the performances of the \basemodel\ as well improve
for increasing $\alpha$, which suggests that the problem is easier when
feasibility constraints are not too demanding. This confirms the outcome of
some observations by \cite{BLM17} about the structure of optimal
solutions of the linear relaxation of this formulation, as these solutions
are allowed to use infeasible paths at a fractional level.

\subsection{Strengthening the model}\label{sec:results-strengthening}

As already mentioned, the \BLM\ approach
is based on a branch-and-cut algorithm in which the arc-flow formulation is iteratively 
strengthened by means of valid inequalities, designed to
cut off infeasible solutions of the linear relaxation.
\cite{BLM17} showed that adding these inequalities
is beneficial to the algorithm, in terms of 
value of the dual bound at the root node and number of instances that 
can be solved to optimality. 

Our third set of experiments is thus aimed at evaluating the impact of 
adding valid inequalities to the path-based formulation.
Table \ref{tab:results-cuts} gives the outcome of our experiments on 
instances in scenarios 1--7 for the \bap\ approach without and with the 
addition of valid inequalities (branch-and-cut-and-price).

\begin{table}[b!]
\centering      
\tabcolsep=4pt
\begin{tabular}{r|rrrr|rrrrrr}
         & \multicolumn{4}{c|}{linear relaxation} & \multicolumn{6}{c}{exact solution} \\
         & \multicolumn{2}{c}{without cuts} & \multicolumn{2}{c|}{with cuts} &                 \multicolumn{3}{c}{\bap} & \multicolumn{3}{c}{\bapC} \\
scenario & \% gap & time & \% gap &   time & \# opt & \% gap &    time &  \# opt  & \% gap &    time \\
\hline
       1 &   5.38 & 0.26 &   3.22 &  31.39 &      5 &   0.00 &   17.90 &        5 &   0.00 &   52.00 \\
       2 &   4.77 & 0.48 &   2.43 & 102.27 &      5 &   0.00 &   63.27 &        5 &   0.00 &  164.91 \\
       3 &   2.79 & 0.49 &   1.24 & 102.82 &      5 &   0.00 &   46.46 &        5 &   0.00 &  146.63 \\
       4 &   6.11 & 0.70 &   3.79 & 186.62 &      5 &   0.00 &  380.83 &        5 &   0.00 &  454.94 \\
       5 &   6.58 & 1.52 &   4.00 & 286.55 &      3 &   1.11 &  220.35 &        3 &   0.98 &  419.83 \\
       6 &   5.02 & 1.60 &   2.64 & 357.34 &      4 &   0.58 &  250.27 &        4 &   0.48 &  549.92 \\
       7 &   7.31 & 2.25 &   4.37 & 600.73 &      4 &   0.99 & 2058.17 &        4 &   0.78 & 1702.62 \\
\hline
 summary &   5.42 & 1.04 &   3.10 & 238.25 &     31 &   0.38 &  401.20 &       31 &   0.32 &  463.29 \\
\end{tabular}
\caption{Results on the addition of valid inequalities.}\label{tab:results-cuts}
\end{table}

The table is organized in two parts. In the first one, we
report the average percentage gap of the linear relaxation in the two
configurations, and the associated computing time reported by SCIP. 
For the version of the algorithm with cuts, we borrowed
from \cite{BLM17} the following families of inequalities: 3OR, 1CUT-IF and 1OR-IF, obtained by combining three OR conditions, one CUT with one or more IF conditions, and one OR with one or more IF conditions, respectively. The reader is referred to \cite{BLM17} for the definition 
of these inequalities as well as to their separation; additional inequalities from this paper showed to have a very marginal effect in our preliminary computational experiments.
Separation is carried out at the root node until no violated cut is found, according
to SCIP tolerance.
The results in Table \ref{tab:results-cuts} confirm that the addition of 
valid inequalities produces a tighter formulation for which the
dual gap with respect to the optimum value is quite small, and reduced by 42\% with
respect to the formulation without cuts (from 5.42\% to 3.10\%).
However, separating these inequalities is time consuming in practice, which
prevents the exhaustive separation of the cuts in an enumerative approach.

For this reason, in the rightmost part of the table we consider 
a \bapC\ algorithm, in which separation is embedded
into the branch-and-price in a heuristic way as follows:
cuts are added at the root node  only, and at most 25 rounds of separation are performed. At each separation round, we consider in order 3OR, 1CUT-IF and 1OR-IF, and we stop the separation as soon as a valid inequality is obtained. 
The inequality is added to the restricted master problem which is then re-optimized.
This heuristic approach is justified by some preliminary experiments on each family of inequalities, where we evaluated the computational effort required for deriving a valid inequality and the relative effect of the inequality on the dual bound.
Remind that a nice property of our approach is that the addition 
of new cuts does not affect the structure of the 
pricing subproblem, yielding a robust branch-and-cut-and-price approach.
 For both \bap\ and 
\bapC\ we report the number of optimal solutions,
the average percentage gap and the average computing time.

The results on the exact methods show that \bapC\ solves the same number of instances
as \bap, and produces slightly better gaps for unsolved instances. Indeed, 
both algorithms solve 31 instances, the average percentage gaps  being 0.38 
(for \bap) and 0.32 (for \bapC).
Despite adding valid cuts seems to be very effective in closing the gap at the root node,
its limited contribution within an enumerative scheme is due to the 
computational overhead required for separating cuts and for solving larger models
at each decision node.

\section{Conclusions}\label{sec:conclusions}

We considered an NP-hard network design problem with end-to-end service 
requirements that play a fundamental role in many contexts, including
telecommunications and transportation.
From a modelling viewpoint, 
we proposed a novel ILP formulation in which 
variables are associated with feasible paths, and discussed alternative
ways for handling the exponential number of variables in the model.
From a methodological perspective, we showed how a column generation 
algorithm can be embedded into a branch-and-cut-and-price scheme, that is 
robust in the sense that  the structure of the subproblems is not 
altered by the branching conditions nor by the addition of valid
inequalities. 
Finally, we gave a comprehensive computational analysis of the 
performances of the proposed algorithm, which is compared with a
state-of-the-art approach proposed in the recent literature. 
Our computational experiments showed that the proposed algorithm
outperforms its competitor and scales efficiently to
larger size of the network.

The introduced path-based formulation is quite general, as all the 
nasty constraints appear in the definition of feasible paths only.
For this reason, it may be worthy to use this modelling approach
for other multi-commodity network design problems arising
in different contexts.

\section{Acknowledgments}

This research was supported by ``Mixed-Integer Non Linear Optimisation: Algorithms and Application'' consortium, which has received funding from the European Union’s EU Framework Programme for Research and Innovation Horizon 2020 under the Marie Sk{\l}odowska-Curie Actions Grant Agreement No. 764759. 
The authors are grateful to Prakash Mirchandani, who kindly provided us the code for generating the numerical data; to Antonio Frangioni, for stimulating discussions;
and to the SCIP optimization suite mailing list,
in particular to Ambros Gleixner, for his invaluable help.

\bibliography{biblio}

\end{sloppypar}
\end{document}